\documentclass{article}

\title{Transient Anomaly Imaging in Visco-Elastic Media Obeying a Frequency
Power-Law}
\author{Elie Bretin\thanks{Centre de
Math{\'e}matiques Appliqu{\'e}es, CNRS UMR 7641, Ecole
Polytechnique, 91128 Palaiseau, France
(bretin@cmap.polytechnique.fr,
lili.guadarrama-bustos@cmap.polytechnique.fr,
wahab@cmap.polytechnique.fr).} \and Lili Guadarrama
Bustos\footnotemark[1] \and Abdul Wahab\footnotemark[1]}

\usepackage[dvips]{graphicx}
\usepackage{graphicx}
\usepackage{amsmath, mathrsfs, wasysym}
\usepackage{amssymb}
\usepackage{textcomp}
\usepackage{stmaryrd}
\usepackage{marvosym}
\usepackage{latexsym, amsmath, amssymb, a4, epsfig}
\usepackage{graphicx}
\usepackage{color}

\newtheorem{thm}{Theorem}[section]

\newtheorem{lem}[thm]{Lemma}

\newcommand{\eqnref}[1]{(\ref {#1})}

\newcommand{\RR}{\mathbb{R}}

\newcommand{\bu}{{\bf{u}}}

\newcommand{\pathfigures}{Figures/}
\graphicspath{{\pathfigures}}

\newcommand{\be}{\begin{equation}}
\newcommand{\ee}{\end{equation}}

\newcommand{\beq}{\begin{equation}}
\newcommand{\eeq}{\end{equation}}

\begin{document}
%%%
%%%
%%%
%%%
\maketitle
%%%
%%%
%%%
%%%
%%%%%%%%%%%%%%%%%%%%
\begin{abstract}%%%%
%%%%%%%%%%%%%%%%%%%%
In this work, we consider the problem of reconstructing a small
anomaly in a viscoelastic medium from wave-field measurements. We
choose Szabo's model
% [Szabo TL, Wu J.  J. Acout. Soc. Am. 2000, 107(5):2437--2446]
to describe the viscoelastic
properties of the medium. Expressing the ideal elastic field
without any viscous effect in terms of the measured field in a
viscous medium, we generalize the imaging procedures, such as time reversal, Kirchhoff Imaging and Back propagation, for an ideal medium  to detect an anomaly in
a visco-elastic medium from wave-field measurements.

\end{abstract}

%%%
%%%
%%%
%%%%%%%%%%%%%%%%%%%%%%%%
\section{Introduction}%%
%%%%%%%%%%%%%%%%%%%%%%%%
%%%
%%%
%%%
We consider the problem of reconstructing a small anomaly in a
viscoelastic medium from wave-field measurements. The Voigt model
is a common model to describe the viscoelastic properties of
tissues. Catheline {\it et al.} \cite{catheline} have shown that
this model is well adapted to describe the viscoelastic response
of tissues to low-frequency excitations. We choose a more general
model derived by Szabo {\it et al.} \cite{SzaboWu00} that
describes observed power-law behavior of many viscoelastic
materials. It is based on a time-domain statement of causality
\cite{Szabo95causal}. It reduces to the Voigt model for the
specific case of quadratic frequency loss. Expressing the ideal
elastic field without any viscous effect in terms of the measured
field in a viscous medium, we  generalize the methods described in
\cite{ HabibEmergBioMedImag,HabibModelBioMed, HabibLiliAcoustic, Habib-Lili10, expansion}; namely the time reversal,
back-propagation and Krichhoff Imaging, to recover the
viscoelastic and geometric properties of an anomaly from wave-field
measurements.

The article is organized as follows. In section \ref{GVEW-Eqn} we
introduce a general visco-elastic wave equation. section \ref{GreenFunc} is devoted to
the derivation of the Green function in a
viscoelastic medium. In section \ref{ImagingProc} we present anomaly
imaging procedures and reconstruction methods in visco-elastic media. Numerical illustrations
are provided in section \ref{Numerics}.
%%%
%%%
%%%
%%%%%%%%%%%%%%%%%%%%%%%%%%%%%%%%%%%%%%%%%%%%%%%%%%%%%%%%%%%%%%%%%%
\section{General Visco-Elastic Wave Equation}\label{GVEW-Eqn} %%%%
%%%%%%%%%%%%%%%%%%%%%%%%%%%%%%%%%%%%%%%%%%%%%%%%%%%%%%%%%%%%%%%%%%
%%%
%%%
%%%
%%%

When a wave travels through a biological medium, its amplitude
decreases with time due to attenuation. The
attenuation coefficient for biological tissue may be approximated
by a power-law over a wide range of frequencies. Measured
attenuation coefficients of soft tissue typically have linear or
greater than linear dependence on frequency \cite{Duck90, Szabo95causal, SzaboWu00}.

In an ideal medium; without attenuation, Hooke's law gives the following relationship between stress and strain tensors:
%%%
%%%
%%%
\begin{equation}
\mathcal{T}=\mathcal{C}:\mathcal{S}
\end{equation}
%%%
%%%
%%%
where $\mathcal{T},~\mathcal{C}$ and $\mathcal{S}$ are respectively stress, stiffness and strain tensors of orders 2, 4 and 2 and : represents tensorial product.

Consider a dissipative medium. Suppose that the medium is homogeneous and isotropic. We write
%%%
%%%
%%%
\begin{eqnarray}
\mathcal{C}=\left[\mathcal{C}_{ijkl}\right]
=\left[\lambda\delta_{ij}\delta_{kl}+\mu(\delta_{ik}\delta_{jl}+\delta_{il}\delta_{jk})\right],
\label{C}
\\
{\eta}=\left[{\eta}_{ijkl}\right] =
\left[\eta_s\delta_{ij}\delta_{kl}+\eta_p(\delta_{ik}\delta_{jl}+\delta_{il}\delta_{jk})\right],
\label{eta}
\end{eqnarray}
%%%
%%%
%%%
where $\delta_{ab}$ is the Kronecker delta function, $\mu,\lambda$ are
the Lam\'e parameters, and $\eta_s,\eta_p$ are the shear and bulk
viscosities, respectively. Here we have adopted the generalized
summation convention over the repeated index.

Throughout this work we suppose that
%%%
%%%
%%%
\begin{equation}
\eta_p, \eta_s <<1. \label{SmallApprox}
\end{equation}
%%%
%%%
%%%

For a medium obeying a power-law attenuation model and under the smallness condition
\eqnref{SmallApprox},  a generalized Hooke's law reads \cite{SzaboWu00}
%%%
%%%
%%%
\begin{equation}
\mathcal{T}(x,t)=\mathcal{C}:\mathcal{S}(x,t)+{\eta}:\mathcal{M}
(\mathcal{S})(x,t) \label{GeneralHookLaw}
\end{equation}
%%%
%%%
%%%
where the convolution operator $\mathcal{M}$ is given by
%%%
%%%
%%%
\begin{equation} \label{defM}
\mathcal{M}(\mathcal{S}) = \left\{
\begin{array}{ll}
-(-1)^{y/2}\frac{\partial^{y-1} \mathcal{S}}{\partial t^{y-1}}
&\text{ y is an even integer},
\\\\
\frac{2}{\pi}(y-1)!(-1)^{(y+1)/2}\frac{H(t)}{t^y}* \mathcal{S}
&\text{ y is an odd integer},
\\\\
-\frac{2}{\pi}\Gamma(y)\sin(y\pi/2)\frac{H(t)}{|t|^y}* \mathcal{S}
&\text{ y is a non integer}.
\end{array}
\right.
\end{equation}
%%%
%%%
%%%
Here $H(t)$ is the Heaviside function and $\Gamma$ denotes the
gamma function.

Note that for the common case, $y=2$, the generalized Hooke's law
(\ref{GeneralHookLaw}) reduces to the Voigt
model,
%%%
%%%
%%%
\begin{equation}
\mathcal{T}=\mathcal{C}:\mathcal{S}+{\eta}:\frac{\partial\mathcal{S}}{\partial
t}. \label{DisHookLaw}
\end{equation}
%%%
%%%
%%%
Taking the divergence of \eqnref{GeneralHookLaw} we get
%%%
%%%
%%%
\begin{equation*}
\nabla\cdot\mathcal{T}
=\left(\bar\lambda+\bar\mu\right)\nabla(\nabla\cdot
\bu)+\bar\mu\Delta \bu,
\end{equation*}
%%%
%%%
%%%
where
%%%
%%%
%%%
\begin{eqnarray*}
\bar\lambda=\lambda+\eta_p\mathcal{M}(\cdot) \quad \mbox{ and }
\quad  \bar\mu=\mu+\eta_s\mathcal{M}(\cdot).
\end{eqnarray*}
%%%
%%%
%%%
Next, considering the equation of motion for the system, {\it
i.e.},
%%%
%%%
%%%
\begin{equation}
\rho\frac{\partial^2 \bu}{\partial t^2}- \mathbf{F}
=\nabla\cdot\mathcal{T}, \label{EqnMotion}
\end{equation}
%%%
%%%
%%%
with $\rho$ being the constant density and $\mathbf{F}$ the
applied force. Using the expression for
$\nabla\cdot\mathcal{T}$, we obtain the generalized visco-elastic
wave equation
%%%
%%%
%%%
\begin{equation} \label{GenVisElasEqn}
\rho\frac{\partial^2 \bu}{\partial t^2}- \mathbf{F}
=\left(\bar\lambda+\bar\mu\right)\nabla(\nabla\cdot
\bu)+\bar\mu\Delta \bu .
\end{equation}
%%%
%%%
%%%
%%%
%%%%%%%%%%%%%%%%%%%%%%%%%%%%%%%%%%%%%%%%%%%%%%%
\section{Green's Function}\label{GreenFunc}%%%%
%%%%%%%%%%%%%%%%%%%%%%%%%%%%%%%%%%%%%%%%%%%%%%%
%%%
%%%
%%%
%%%
In this section we find the Green function of the viscoelastic
wave equation \eqnref{GenVisElasEqn}. For doing so, we first need
a Helmholtz decomposition.

%%%
%%%
%%%
%%%
%%%%%%%%%%%%%%%%%%%%%%%%%%%%%%%%%%%%%%
\subsection{Helmholtz Decomposition}%%
%%%%%%%%%%%%%%%%%%%%%%%%%%%%%%%%%%%%%%
%%%
%%%
%%%
The following lemma holds.
%%%
%%%
%%%%%%%%%%%%%%%%%%%%%%%%%%
\begin{lem}\label{LAME}%%%
%%%%%%%%%%%%%%%%%%%%%%%%%%
%%%
%%%
If the displacement field $\bu(x,t)$ satisfies
\eqnref{GenVisElasEqn}, $\frac{\partial\bu(x,0)}{\partial t}=\nabla A+\nabla\times B$ and $\bu(x,0)=\nabla C+\nabla\times D$  and if the body force
$\mathbf{F}=\nabla\varphi_f+\nabla\times{\psi}_f$ then
there exist potentials $\varphi_u$  and ${\psi}_u$ such that
%%%
%%%
%%%
%%%
\begin{itemize}
%%%
%%%
%%%
\item $\bu =\nabla\varphi_u+\nabla\times{\psi}_u$;
$\nabla\cdot{\psi}_u=0$;
%%%
%%%
%%%
\item
$\frac{\partial^2\varphi_u}{\partial t^2}
=\frac{\varphi_f}{\rho}+c_p^2\Delta \varphi_u+
\nu_p\mathcal{M}(\Delta \varphi_u) \approx \frac{\varphi_f}{\rho}
- \frac{\nu_p \mathcal{M}(\varphi_f)}{\rho c_p^2}
+c_p^2\Delta \varphi_u+ \frac{\nu_p}{c_p^2} \mathcal{M}(
\partial_t^2\varphi_u)$;
%%%
%%%
%%%
\item $\frac{\partial^2{\psi}_u}{\partial t^2}
=\frac{{
\psi}_f}{\rho}+c_s^2\Delta{\psi}_u+\nu_s\mathcal{M}(\Delta
{\psi}_u)\approx \frac{{ \psi}_f}{\rho} -
\frac{\nu_s \mathcal{M}({ \psi}_f)}{\rho c_s^2}
+c_s^2\Delta{\psi}_u+\frac{\nu_s}{c_s^2} \mathcal{M}(
\partial_t^2{\psi}_u)$,
\end{itemize}
%%%
%%%
%%%
with $$c_p^2=\frac{\lambda+2\mu}{\rho},~c^2_s=\frac{\mu}{\rho},
~\nu_p=\frac{\eta_p+2\eta_s}{\rho}, \quad \mbox{and } \quad
\nu_s=\frac{\eta_s}{\rho}.$$
%%%
%%%
%%%%%%%%%%%%
\end{lem}%%%
%%%%%%%%%%%%
%%%
%%%
%%%%%%%%%%%%%%%%%%%%%%%%%%
\textit{\textbf{Proof.}}%%
%%%%%%%%%%%%%%%%%%%%%%%%%%
%%%
%%%
For $\varphi_u$ and $\psi_u$ defined as
%%%
%%%
%%%
\begin{eqnarray}
\varphi_u(x,t)
=\int_0^t\int_0^\tau
\left[\frac{\varphi_f}{\rho}+(c_p^2+\nu_p\mathcal{M})(\nabla\cdot u)\right]dsd\tau+tA+C\label{Phi_u}
\\
\psi_u(x,t)
=\int_0^t\int_0^\tau\left[\frac{\vec\psi_f}{\rho}-(c_s^2+\nu_s\mathcal{M})(\nabla\times u)\right]dsd\tau
+t\vec B+\vec D\label{Psi_u}
\end{eqnarray}
%%%
%%%
%%%
%%%
we have the required expression for $\bu$. Moreover, it is evident from (\ref{Psi_u}) that $\nabla\cdot\psi_u=0$\\

Now, on differentiating $\varphi_u$ and $\psi_u$ twice with respect to time, we get
%%%
%%%
%%%
\begin{equation*}
\frac{\partial^2\varphi_u}{\partial t^2}
=\frac{\varphi_f}{\rho}+c_p^2\Delta\varphi_u+ \nu_p\mathcal{M}(\Delta\varphi_u)
\end{equation*}
%%%
%%%
%%%
\begin{equation*}
\frac{\partial^2{\psi}_u}{\partial t^2}
=\frac{{\psi}_f}{\rho}+c_s^2\Delta{\psi}_u+\nu_s\mathcal{M}(\Delta
{\psi}_u)
\end{equation*}
%%%
%%%
%%%
Finally, applying $\mathcal{M}$ on last two equations, neglecting the higher order terms in $\nu_s$ and $\nu_p$ and injecting back the expressions for $\mathcal{M}(\Delta\varphi_u)$ and $\mathcal{M}(\Delta\psi_u)$, we get the required differential equations for $\varphi_u$ and $\psi_u$. \hfil$\Box$\\

Let
%%%
%%%
%%%
\be \label{k2M} K_{m}(\omega) = \omega \sqrt{ \left(1-
\frac{\nu_{m}}{c_{m}^2} \hat{\mathcal{M}}(\omega)\right)}, \quad
m=s,p,\ee
%%%
%%%
%%%
 where the multiplication operator
$\hat{\mathcal{M}}(\omega)$ is the Fourier transform of the
convolution operator $\mathcal{M}$.

If $\varphi_u$ and ${\psi}_u$ are causal then it implies the causality of the
inverse Fourier transform of $K_{m}(\omega), m=s,p$. Applying the
Kramers-Kr\"{o}nig relations, it follows that
%%%
%%%
%%%
\be \label{KKr} - \Im m K_m(\omega) =
\mathcal{H}\bigg[ \Re e K_m(\omega) \bigg] \quad \mbox{and} \quad
\Re e  K_m(\omega) = \mathcal{H}\bigg[ \Im m K_m(\omega) \bigg],
\quad m=p,s, \ee
%%%
%%%
%%%
where $\mathcal{H}$ is the Hilbert transform.
Note that $\mathcal{H}^2 = -I$. The convolution operator
$\mathcal{M}$ given by \eqnref{defM} is based on the constraint
that causality imposes on \eqnref{GeneralHookLaw}. Under the
smallness assumption \eqnref{SmallApprox}, the expressions in
\eqnref{defM} can be found from the Kramers-Kr\"{o}nig relations
\eqnref{KKr}. One drawback of \eqnref{KKr} is that the
\index{attenuation} attenuation, $\Im m K_m(\omega)$, must be
known at all frequencies to determine \index{dispersion} the
dispersion, $\Re e K_m(\omega)$. However, bounds on the dispersion
can be obtained from measurements of the attenuation over a finite
frequency range \cite{milton}.
%%%
%%%
%%%
%%%

%%%%%%%%%%%%%%%%%%%%%%%%%%%%%%%%%%%%%%%%%%%
\subsection{Solution of (\ref{GenVisElasEqn}) with a Concentrated Force.}%%%
%%%%%%%%%%%%%%%%%%%%%%%%%%%%%%%%%%%%%%%%%%%
%%%
%%%
%%%
Let $u_{ij}$ denote the $i$-th component of the solution $\bu_j$
of the elastic wave equation related to a force $\mathbf{F}$
concentrated in the $x_j$-direction. Let $j=1$ for simplicity and
suppose that
%%%
%%%
%%%
\begin{equation}
\mathbf{F} =-T(t)\delta(x-\xi)\mathbf{e}_1
=-T(t)\delta(x-\xi)(1,0,0), \label{ConsentF}
\end{equation}
%%%
%%%
%%%
where $\xi$ is the source point and
$(\mathbf{e}_1,\mathbf{e}_2,\mathbf{e}_3)$ is an orthonormal basis
of $\RR^3$. The corresponding Helmholtz decomposition of the force
$\mathbf{F}$ can be written \cite{ElasticWavePujol} as

%%%
%%%
%%%
\begin{equation}
\left\{
\begin{array}{l}
\mathbf{F}=\nabla\varphi_f+\nabla\times{\psi}_f,
\\
\\
\varphi_f=\frac{T(t)}{4\pi}\frac{\partial}{\partial
x_1}\left(\frac{1}{r}\right),
\\
\\
{\psi}_f=-\frac{T(t)}{4\pi} \left( 0
,\frac{\partial}{\partial x_3}\left(\frac{1}{r}\right)
,-\frac{\partial}{\partial x_2}\left(\frac{1}{r}\right) \right),
\end{array}
\right.\label{HelConsentF}
\end{equation}
%%%
%%%
%%%
where $r=|x-\xi|$.

Consider the Helmholtz decomposition for $\bu_{i1}$ as
%%%
%%%
%%%
\begin{equation}
u_{i1}= \nabla\varphi_1+\nabla\times\vec\psi_1\label{ui1}
\end{equation}
%%%
%%%
%%%
where $\varphi_1$ and $\psi_1$ are the solutions of the equations
%%%
%%%
%%%
\begin{eqnarray}
\Delta
\varphi_1-\frac{1}{c^2_p}\frac{\partial^2\varphi_1}{\partial t^2}
+ \frac{\nu_p}{c_p^4} \mathcal{M}(\partial_t^2 \varphi_1)
=\frac{\nu_p \mathcal{M}({
\varphi}_f)}{\rho c_p^4}-\frac{\varphi_f}{c_p^2\rho}, \label{PhiDiffCon}
\\
\Delta
{\psi}_1-\frac{1}{c^2_s}\frac{\partial^2{\psi}_1}{\partial
t^2} + \frac{\nu_s}{c_s^4}\mathcal{M}(\partial_t^2{\psi}_1)
=\frac{\nu_s\mathcal{M}({ \psi}_f)}{\rho c_s^4}-\frac{{\psi}_f}{c_s^2\rho}. \label{PsiDiffCon}
\end{eqnarray}
%%%
%%%
%%%

Taking the Fourier transform of (\ref{ui1}),(\ref{PhiDiffCon}) and
(\ref{PsiDiffCon}) with respect to $t$ we get
%%%
%%%
%%%
\begin{eqnarray}
\hat \bu_{1}
=\nabla\hat\varphi_1+\nabla\times\hat\psi_1\label{ui1Hat}
\\
\Delta \hat\varphi_1+\frac{K^2_p(\omega)}{c_p^2}\hat\varphi_1
=\frac{\nu_p
\hat{\mathcal{M}}(\omega) \hat{{ \varphi}}_f}{\rho c_p^4}- \frac{\hat\varphi_f}{\rho c_p^2},
\label{phi1hel}
\\
\Delta \hat\psi_1+\frac{K^2_s(\omega)}{c_s^2}\hat\psi_1
=\frac{\nu_s
\hat{\mathcal{M}}(\omega)\hat{\psi}_f}{\rho c_s^4}-\frac{\hat\psi_f}{\rho c_s^2},
\label{psi1hel}
\end{eqnarray}
%%%
%%%
%%%
with $ K_m(\omega), m=p,s,$ given by \eqnref{k2M}.

It is well known that the Green's functions of the Helmholtz equations (\ref{phi1hel}) and (\ref{psi1hel}) are
%%%
%%%
%%%
\begin{eqnarray*}
\hat g^m(r,\omega)
=\frac{e^{\sqrt{-1}\frac{K_m(\omega)}{c_m}r}}{4\pi r}, \quad
m=s,p.
\end{eqnarray*}
%%%
%%%
%%%
Thus, following \cite{ElasticWavePujol} we write $\hat\varphi_1$ as
%%%
%%%
%%%
\begin{equation*}
\hat\varphi_1(x,\omega;\xi) =-\left(1-
\frac{\nu_p \hat{\mathcal{M}}(\omega) }{c_p^2}\right) \frac{\hat
T(\omega)}{\rho(4\pi c_p)^2}\int_V\hat g^p(x-\chi,\omega)\frac{\partial}{\partial\chi_1}\frac{1}{|\chi-\xi|}dV_\chi.
\end{equation*}
%%%
%%%
%%%
and divide the volume $V$ into spherical shells of radius $h$ centered at observation point $x$. On each shell $\hat g^p(x-\chi,\omega)$ rests constant. So we have
%%%
%%%
%%%
\begin{equation*}
\hat\varphi_1(x,\omega;\xi) =-\left(1-
\frac{\nu_p \hat{\mathcal{M}}(\omega) }{c_p^2}\right) \frac{\hat
T(\omega)}{\rho(4\pi c_p)^2}\int_0^\infty\frac{1}{h}\hat g^p(h,\omega)\int_\sigma\frac{\partial}{\partial\chi_1}\frac{1}{R}d\sigma dh.
\end{equation*}
%%%
%%%
%%%
with $h=|x-\chi|$, $R=|\chi-\xi|$ and $d\sigma$ the appropriate surface element.

As \cite{AkiRichard}
%%%
%%%
%%%
\begin{equation*}
\int_\sigma\frac{\partial}{\partial\chi_1}\left(\frac{1}{R}\right)
=\left\{
\begin{array}{ll}
0 &h>r
\\
4\pi h^2\frac{\partial}{\partial x_1}\left(\frac{1}{r}\right) &h<r
\end{array}
\right.\label{SurfaceInt}
\end{equation*}
%%%
%%%
%%%
Therefore, we have following expression for  $\hat\varphi_1$:
%%%
%%%
%%%
\begin{equation}
\hat\varphi_1(x,\omega;\xi) =-\left(1-
\frac{\nu_p \hat{\mathcal{M}}(\omega) }{c_p^2}\right) \frac{\hat
T(\omega)}{4\pi\rho}\frac{\partial}{\partial
x_1}\left(\frac{1}{r}\right) \int_0^{r/c_p} \zeta e^{\sqrt{-1}
K_p(\omega)\zeta}\, d\zeta . \label{phiHatValue}
\end{equation}
%%%
%%%
%%%
In the same way, the  vector $\hat\psi_1$ is given by
%%%
%%%
%%%
\begin{equation}
\hat\psi_1(x,\omega;\xi) =\left(1- \frac{\nu_s
\hat{\mathcal{M}}(\omega) }{c_s^2}\right)\frac{\hat T(\omega)}{4\pi\rho}
\left( 0 ,\frac{\partial}{\partial x_3}\left(\frac{1}{r}\right)
,-\frac{\partial}{\partial x_2}\left(\frac{1}{r}\right) \right)
\int_0^{r/c_s}\zeta e^{\sqrt{-1} K_s(\omega)\zeta}\, d\zeta.
\label{psiHatValue}
\end{equation}
%%%
%%%
%%%
Introduce the following notation:
%%%
%%%
%%%
\begin{eqnarray}
I_m(x,\omega)=A_m \int_0^{r/c_m} \zeta e^{\sqrt{-1}
K_m(\omega)\zeta}\, d\zeta
\\
E_m(x,\omega)= A_m e^{\sqrt{-1} K_m(\omega)\frac{r}{c_m}}, \\
A_m(\omega) =  \left(1- \frac{\nu_m \hat{\mathcal{M}}(\omega)
}{c_m^2}\right), \quad m=p,s.
\end{eqnarray}
%%%
%%%
%%%
We obtain, after a lengthy but simple calculation, that $\hat
u_{i1}$ is given by
%%%
%%%
%%%
\begin{equation*}
\begin{array}{l}
\hat u_{i1} = \frac{\hat
T(\omega)}{4\pi\rho}\frac{\partial^2}{\partial
x_ix_1}\left(\frac{1}{r}\right)\left[I_s(r,\omega)-I_p(r,\omega)\right]
+\frac{\hat T(\omega)}{4\pi \rho c_p^2 r}\frac{\partial
r}{\partial x_i}\frac{\partial r}{\partial x_1}E_p(r,\omega)
\\
\\
\qquad \qquad +\frac{\hat T(\omega)}{4\pi\rho
c_s^2r}\left(\delta_{i1}-\frac{\partial r}{\partial
x_i}\frac{\partial r}{\partial x_1}\right) E_s(r,\omega),
\end{array}
\end{equation*}
%%%
%%%
%%%
and therefore, it follows that the solution $u_{ij}$ for an
arbitrary $j$ is
%%%
%%%
%%%
\begin{equation*}
\begin{array}{l}
\hat u_{ij} =\frac{\hat
T(\omega)}{4\pi\rho}\left(3\gamma_i\gamma_j-\delta_{ij}\right)\frac{1}{r^3}\left[I_s(r,\omega)-I_p(r,\omega)\right]
+\frac{\hat T(\omega)}{4\pi\rho
c_p^2}\gamma_i\gamma_j\frac{1}{r}E_p(r,\omega)
\\
\\
\qquad \qquad +\frac{\hat T(\omega)}{4\pi\rho
c_s^2}\left(\delta_{ij}-\gamma_i\gamma_j\right)\frac{1}{r}E_s(r,\omega),
\end{array}
\end{equation*}
%%%
%%%
%%%
where $\gamma_i = (x_i-\xi_i)/r$.

%%%
%%%
%%%
%%%
%%%%%%%%%%%%%%%%%%%%%%%%%%%%%%%%%
\subsection{Green's function }%%%
%%%%%%%%%%%%%%%%%%%%%%%%%%%%%%%%%
%%%
%%%
%%%
If we substitute $T(t)=\delta(t)$, where delta is the Dirac mass,
then the function $u_{ij}=G_{ij}$ is the $i$-th component of the
Green function related to the force concentrated in the
$x_j$-direction. In this case, we have $\hat T(\omega)=1$. Thus, we
have the following expression for $\hat G_{ij}$:
%%%
%%%
%%%
\begin{equation*}
\begin{array}{l}
\hat
G_{ij}=\frac{1}{4\pi\rho}\left(3\gamma_i\gamma_j-\delta_{ij}\right)\frac{1}{r^3}\left[I_s(r,\omega)-I_p(r,\omega)\right]
+\frac{1}{4\pi\rho c_p^2}\gamma_i\gamma_j\frac{1}{r}E_p(r,\omega)
\\
\\
\qquad \qquad +\frac{1}{4\pi\rho
c_s^2}\left(\delta_{ij}-\gamma_i\gamma_j\right)\frac{1}{r}E_s(r,\omega),
\end{array}
\end{equation*}
%%%
%%%
%%%
which implies that
%%%
%%%
%%%
\be \hat G_{ij}(r,\omega;\xi)=\hat
g^p_{ij}(r,\omega)+\hat g^s_{ij}(r,\omega)+\hat
g^{ps}_{ij}(r,\omega), \label{GreenFonVEWE-Fourier} \ee
%%%
%%%
%%%
 where
%%%
%%%
%%%
\be
 \hat
g^{ps}_{ij}(r,\omega)
=\frac{1}{4\pi\rho}\left(3\gamma_i\gamma_j-\delta_{ij}\right)
\frac{1}{r^3}\left[I_s(r,\omega)-I_p(r,\omega)\right],\label{gps_ij-fourier}
\ee
%%%
%%%
%%%
\be \hat g^p_{ij}(r,\omega) =\frac{A_p(\omega)}{\rho
c_p^2}\gamma_i\gamma_j\hat g^p(r,\omega), \label{gp_ij-Fourier}
\ee and \be \hat g^s_{ij}(r,\omega) =\frac{A_s(\omega)}{\rho
c_s^2}\left(\delta_{ij}-\gamma_i\gamma_j\right)\hat
g^s(r,\omega).\label{gs_ij-Fourier} \ee
%%%
%%%
%%%

Let $G(r,t;\xi)=\left(G_{ij}(r,t;\xi)\right)$ denote the transient
Green function of \eqnref{GenVisElasEqn} associated with the
source point $\xi$. Let $G^m(r,t;\xi)$ and $W_m(r,t)$ be the
inverse Fourier transforms of $A_m(\omega)\hat{g}^m(r,\omega)$ and
$I_m(r,\omega), m=p,s$, respectively. Then, from
(\ref{GreenFonVEWE-Fourier}-\ref{gs_ij-Fourier}), we have
%%ù
%%%
%%%
\begin{equation} \label{Green_function}
\begin{array}{l}
G_{ij}(r,t;\xi) =\frac{1}{\rho c_p^2}\gamma_i\gamma_jG^p(r,t;\xi)
+\frac{1}{\rho
c_s^2}\left(\delta_{ij}-\gamma_i\gamma_j\right)G^s(r,t;\xi)
\\
\\
\qquad \qquad
+\frac{1}{4\pi\rho}\left(3\gamma_i\gamma_j-\delta_{ij}\right)
\frac{1}{r^3}\left[W_s(r,t)-W_p(r,t)\right].
\end{array}
\end{equation}
%%%
%%%
%%%
Note that by a change of variables,
%ù%
%%%
%%%
$$
W_m(r,t)
=\frac{4\pi}{c_m^2}\int_0^{r}\zeta^2G^m(\zeta,t;\xi)d\zeta.
$$
%%%
%%%
%%%
%%%%%%%%%%%%%%%%%%%%%%%%%%%%%%%%%%%%%%%%%%%%%%%%%%
\section{Imaging procedure} \label{ImagingProc}%%%
%%%%%%%%%%%%%%%%%%%%%%%%%%%%%%%%%%%%%%%%%%%%%%%%%%
%%%
%%%
%%%
Consider the limiting case $\lambda \rightarrow +\infty$. The
Green function for a quasi-incompressible visco-elastic medium is
given by
%%%
%%%
%%%
\begin{equation*}
\begin{array}{l}
G_{ij}(r,t;\xi) = \frac{1}{\rho
c_s^2}\left(\delta_{ij}-\gamma_i\gamma_j\right)G^s(r,t;\xi)
\\
\\
\qquad \qquad +\frac{1}{\rho
c_s^2}\left(3\gamma_i\gamma_j-\delta_{ij}\right) \frac{1}{r^3}
\int_0^{r}\zeta^2G^s(\zeta,t;\xi)d\zeta.
\end{array}
\end{equation*}
%%%
%%%
%%%
To generalize the detection algorithms presented in \cite{HabibEmergBioMedImag, HabibModelBioMed, Habib-Lili10, HabibLiliAcoustic,expansion} to the visco-elastic case we shall express the ideal
Green function without any viscous effect in terms of the Green
function in a viscous medium. From
%%%
%%%
%%%
$$
G^s(r,t;\xi) = \frac{1}{\sqrt{2\pi}} \int_\RR e^{- \sqrt{-1}
\omega t} A_s(\omega) g^s(r,\omega)\, d\omega,
$$
%%%
%%%
%%%
it follows that
%%%
%%%
%%%
$$
 G^s(r,t;\xi) = \frac{1}{\sqrt{2\pi}} \int_\RR A_s(\omega)
\frac{e^{\sqrt{-1} (- \omega t+ \frac{K_s(\omega)}{c_s}r)}} {4\pi
r} \, d\omega.
$$
%%%
%%%
%%%
%%%%%%%%%%%%%%%%%%%%%%%%%%%%%%%%%%%%%%%%%%%%%%%%%%
\subsection{Approximation of the Green Function}%%
%%%%%%%%%%%%%%%%%%%%%%%%%%%%%%%%%%%%%%%%%%%%%%%%%%
%%%
%%%
%%%
Introduce the operator
$$L \phi(t) = \frac{1}{ 2\pi }  \int_\RR \int_0^{+\infty} A_s(\omega) \phi(\tau) e^{\sqrt{-1}
K_s(\omega) \tau} e^{- \sqrt{-1} \omega t}\, d\tau \, d\omega,$$
for a causal function $\phi$. We have
%%%
%%%
%%%
$$
G^s(r,t;\xi) = L(\frac{\delta(\tau - r/c_s)}{4\pi r}),
$$
%%%
%%%
%%%
and therefore,
%%%
%%%
%%
$$
L^*G^s(r,t;\xi) = L^*L(\frac{\delta(\tau - r/c_s)}{4\pi r}),
$$
%%ù
%ù%
%%%
where $L^*$ is the $L^2(0,+\infty)$-adjoint of $L$.

Consider for simplicity \index{Voigt model} the Voigt model. Then,
$\hat{\mathcal{M}}(\omega) = - \sqrt{-1} \omega$ and hence,

%%%
%%%
%%%
$$K_s(\omega) = \omega \sqrt{1+ \frac{\sqrt{-1} \nu_s}{c_s^2}\omega}
\approx \omega + \frac{\sqrt{-1} \nu_s}{2 c_s^2} \omega^2,$$
%%%
%%%
%%%
under the smallness assumption \eqnref{SmallApprox}. The operator $L$
can then be approximated by
%%%
%%%
%%%
$$ \tilde L \phi(t) = \frac{1}{ 2\pi }  \int_\RR \int_0^{+\infty} A_s(\omega) \phi(\tau)
e^{- \frac{\nu_s}{2 c_s^2} \omega^2
 \tau} e^{\sqrt{-1} \omega (\tau -t)}\, d\tau \, d\omega.$$
%%%
%%%
%%%
Since

%%%
%%%
%%%
$$
 \int_\RR e^{- \frac{\nu_s}{2 c_s^2} \omega^2
 \tau} e^{\sqrt{-1} \omega (\tau -t)}\,  d\omega =
 \frac{\sqrt{2\pi} c_s}{\sqrt{\nu_s \tau}} e^{-
 \frac{c_s^2(\tau-t)^2}{2\nu_s \tau}},
$$
%%%
%%%
%%%
and
%%%
%%%
%%%
$$
\sqrt{-1} \int_\RR  \omega e^{- \frac{\nu_s}{2 c_s^2} \omega^2
 \tau} e^{\sqrt{-1} \omega (\tau -t)}\,  d\omega = -
 \frac{\sqrt{2\pi} c_s}{\sqrt{\nu_s \tau}} \frac{\partial}{\partial t} e^{-
 \frac{c_s^2(\tau-t)^2}{2\nu_s \tau}},
 $$
%%%
%%%
%%%
it follows that
%%%
%%%
%%%
\be \tilde L \phi(t) = \int_0^{+\infty}  \frac{t}{\tau} \phi(\tau)  \frac{ c_s}{\sqrt{2 \pi \nu_s \tau}} e^{-
 \frac{c_s^2(\tau-t)^2}{2\nu_s \tau}} \, d\tau.\ee
%%%
%%%
%%ù
 Analogously,

%%%
%%%
%%%
\be \label{phase}  \tilde L^* \phi(t) =  \int_0^{+\infty} \frac{\tau}{t}\phi(\tau)
 \frac{c_s}{\sqrt{2 \pi \nu_s t}} e^{-
 \frac{c_s^2(\tau-t)^2}{2\nu_s t}}
\, d\tau.\ee
%%%
%%%
%%%

Since the phase in \eqnref{phase} is quadratic and
$\nu_s$ is small then by the stationary phase theorem \ref{StationaryPhase}, we can prove the following theorem:
\begin{thm}\label{Thm:Approx}
%%%
%%%
%%%
$$ \tilde L^*\phi \approx \phi +  \frac{\nu_s}{2 c_s^2}  \partial_{tt}  ( t \phi  ), \quad \tilde L \phi \approx \phi +  \frac{\nu_s}{2 c_s^2}  t \partial_{tt} \phi,    $$
%%%
%%%
%%%
and
%%%
%%%
%%%
\begin{equation} \label{formasympfinal} \tilde L^* \tilde L\phi
\approx \phi + \frac{\nu_s}{c_s^2}
\partial_t (t
\partial_t \phi), \end{equation}
%%%
%%%
%%%
and therefore,
\begin{equation} \label{formasympfinalinverse}
(L^* \tilde L)^{-1} \phi \approx \phi - \frac{\nu_s}{c_s^2}
\partial_t (t
\partial_t \phi).\end{equation} $\Box$
\end{thm}
{\bf{Proof.}} (See appendix \ref{ApenA})
%%%%%%%%%%%%%%%%%%%%%%%%%%%%%%%%%%%%%%%%%%%%%%%%%%%%%%%%%%%%%%%
\subsection{Reconstruction Methods}\label{ReconstrucMethods}
%%%%%%%%%%%%%%%%%%%%%%%%%%%%%%%%%%%%%%%%%%%%%%%%%%%%%%%%%%%%%%%
%%%
%%%

From the previous section, it follows that the ideal Green
function, ${\delta(\tau -r/c_s)}/{(4\pi r)}$, can be approximately
reconstructed from the viscous Green function, $G^s(r,t;\xi)$, by
either solving the ODE
%%%
%%%
%%%
$$
\phi + \frac{\nu_s}{c_s^2} \partial_t (t
\partial_t \phi) = L^*G^s(r,t;\xi),
$$
%%%
%%%
%%%
with $\phi=0, t \ll 0$ or just making the approximation
$$
{\delta(\tau -r/c_s)}/{(4\pi r)} \approx L^*G^s(r,t;\xi) -
\frac{\nu_s}{c_s^2}
\partial_t (t
\partial_t L^*G^s(r,t;\xi)).
$$

Once the ideal Green function ${\delta(\tau -r/c_s)}/{(4\pi r)}$
is reconstructed, one can find its source $\xi$ using a
time-reversal, a Kirchhoff or a back-propagation algorithm. See
\cite{HabibEmergBioMedImag,
 HabibModelBioMed,HabibLiliAcoustic, Habib-Lili10}.

Using the asymptotic formalism developed in \cite{Habib-Lili10, bookAK,
AKbook}, one can also find the shear modulus of the
anomaly using the ideal near-field measurements which can be
reconstructed from the near-field measurements in the viscous
medium. The asymptotic formalism reduces the anomaly imaging
problem to the detection of the location and the reconstruction of
a certain polarizability tensor in the far-field  and separates
the scales in the near-field.
%%%
%%%
%%%
%%%%%%%%%%%%%%%%%%%%%%%%%%%%%%
\section{Numerical Illustrations}\label{Numerics}
%%%%%%%%%%%%%%%%%%%%%%%%%%%%%%ù
%%%
%%%
%%%
In this section, we illustrate the profile of the Green function. We choose parameters of simulation as in the work of Bercoff {\it et al.}\cite{Bercoff04}: we take $\rho = 1000$, $c_s = 1$, $c_p = 40$, $r = 0.015$ and $\nu_p = 0$. \\

In figure 1, we plot temporal representation of the green function: 
$$ t \to \frac{1}{\rho c_p^2} \left( G^p(r,t;\xi) + G^s(r,t;\xi) \right) +  \frac{1}{4\pi\rho r^3}\left[W_s(r,t)-W_p(r,t)\right].$$  
for three different values of $y$ and $\nu_s$. We can see that the attenuation behavior varies with respect to different choices of power law exponent $y$. One can clearly distinguish the three different terms of the Green function; \emph{i.e.} $G^s_{ij},~G^p_{ij}$ and $G^{ps}_{ij}$.

Figure 2 corresponds to spatial representation of the  green function: 
$$ (x,y) \to \frac{1}{\rho c_p^2} \left( (x/r)^2 G^p(r,t;\xi) + (1-(x/r)^2) G^s(r,t;\xi) \right) +  \frac{1}{4\pi\rho r^3} (3 (x/r)^2-1 )\left[W_s(r,t)-W_p(r,t)\right], $$ 
for different values of $y$ at $t=0.015$. As expected, we get a diffusion of the wavefront with the increasing values of $y$ and depending the choice of $\nu_s$.

In figure 3, we illustrate the results of the approximation of the operator $L\phi$ with the smooth function $\phi(t) = exp(-50*(t-1).^2)''$. As shown by the stationary phase theorem \ref{StationaryPhase} , the numerically calculated $L^{\infty}$-error 
$$ \|L \phi - \left( \phi + \frac{\nu_s}{2 c_s^2} t \phi''\right)  \|_{L^{\infty}(\RR^+)}$$ is of order two.
  
 \begin{figure}  
\begin{center} \label{fig:green_t}
  \includegraphics[width=6cm]{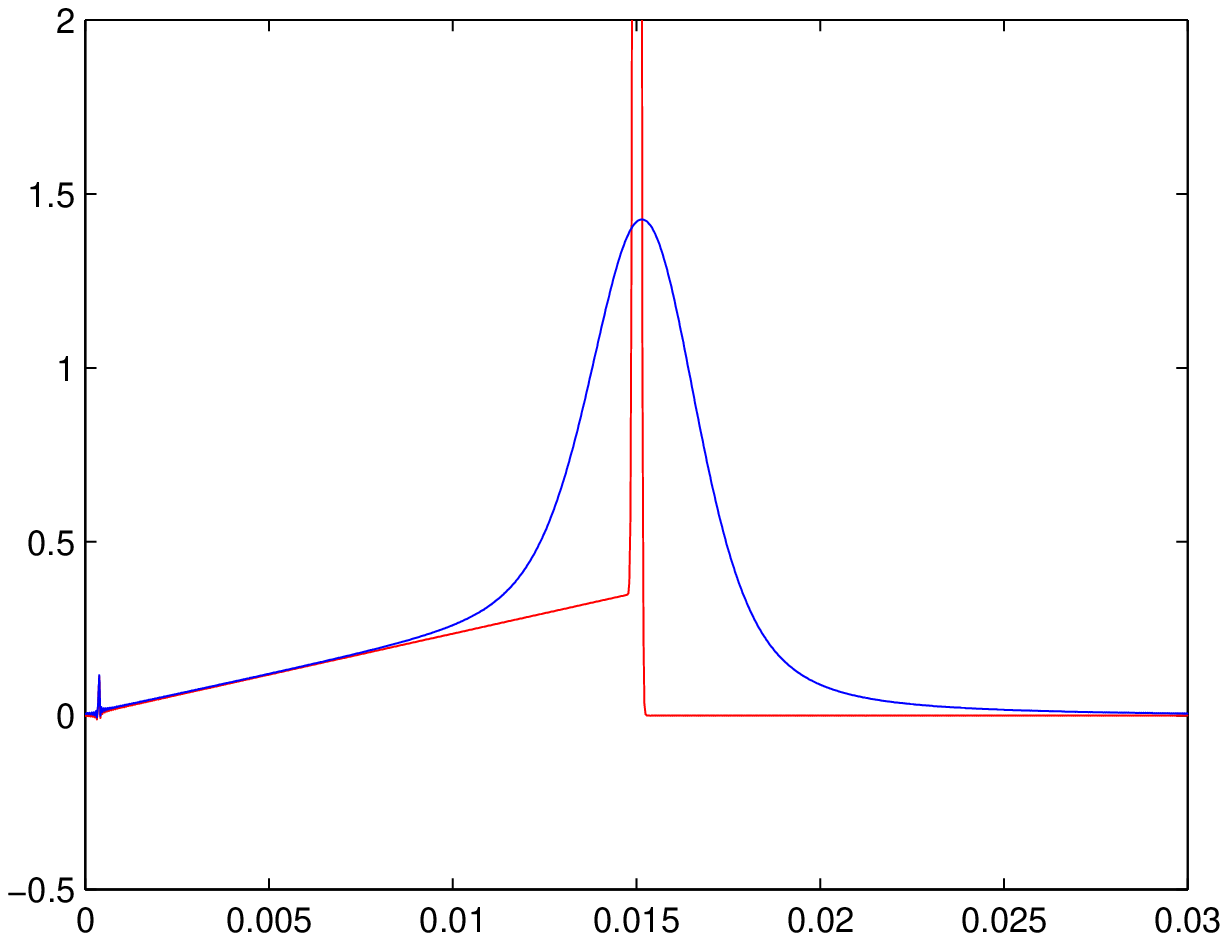}
 \includegraphics[width=6cm]{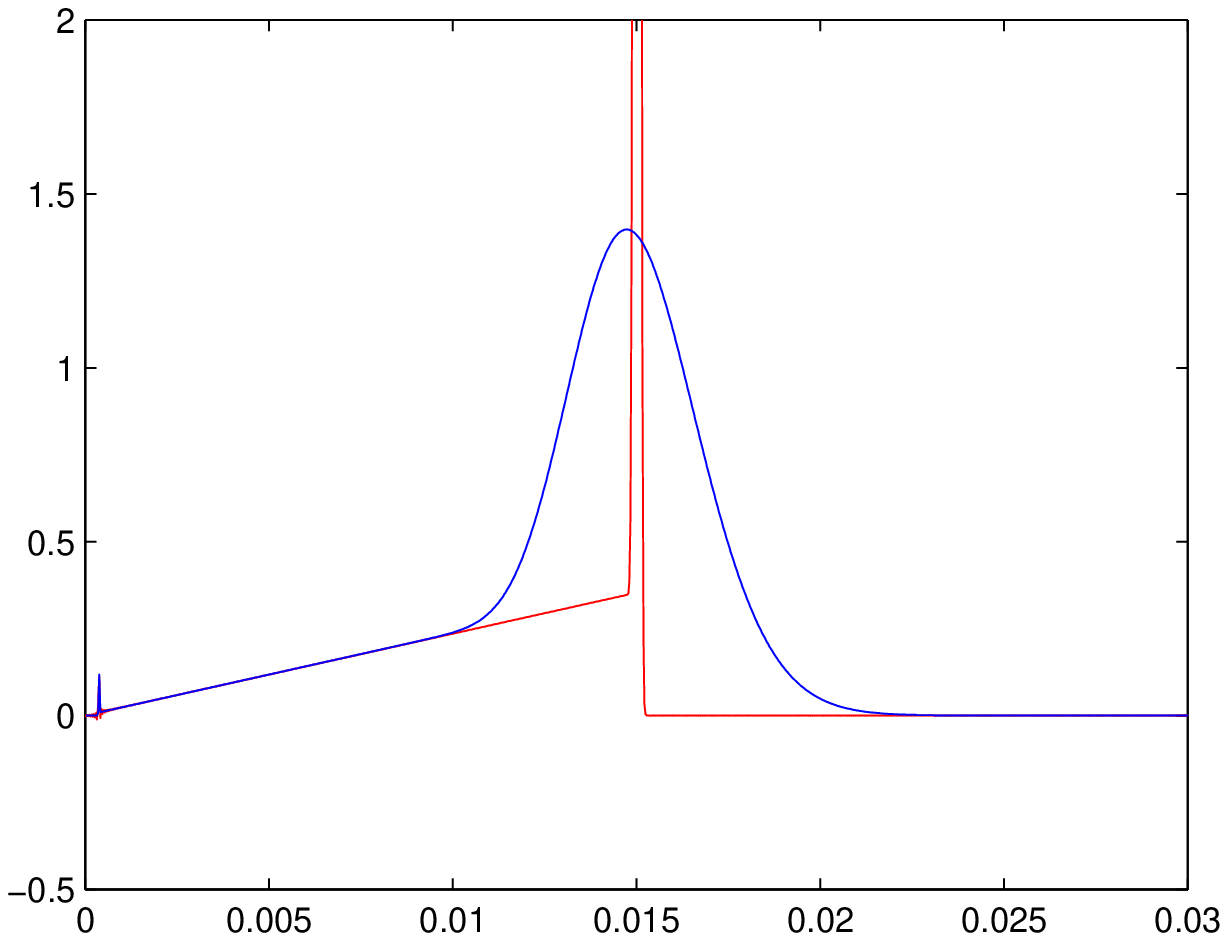}
 \includegraphics[width=6cm]{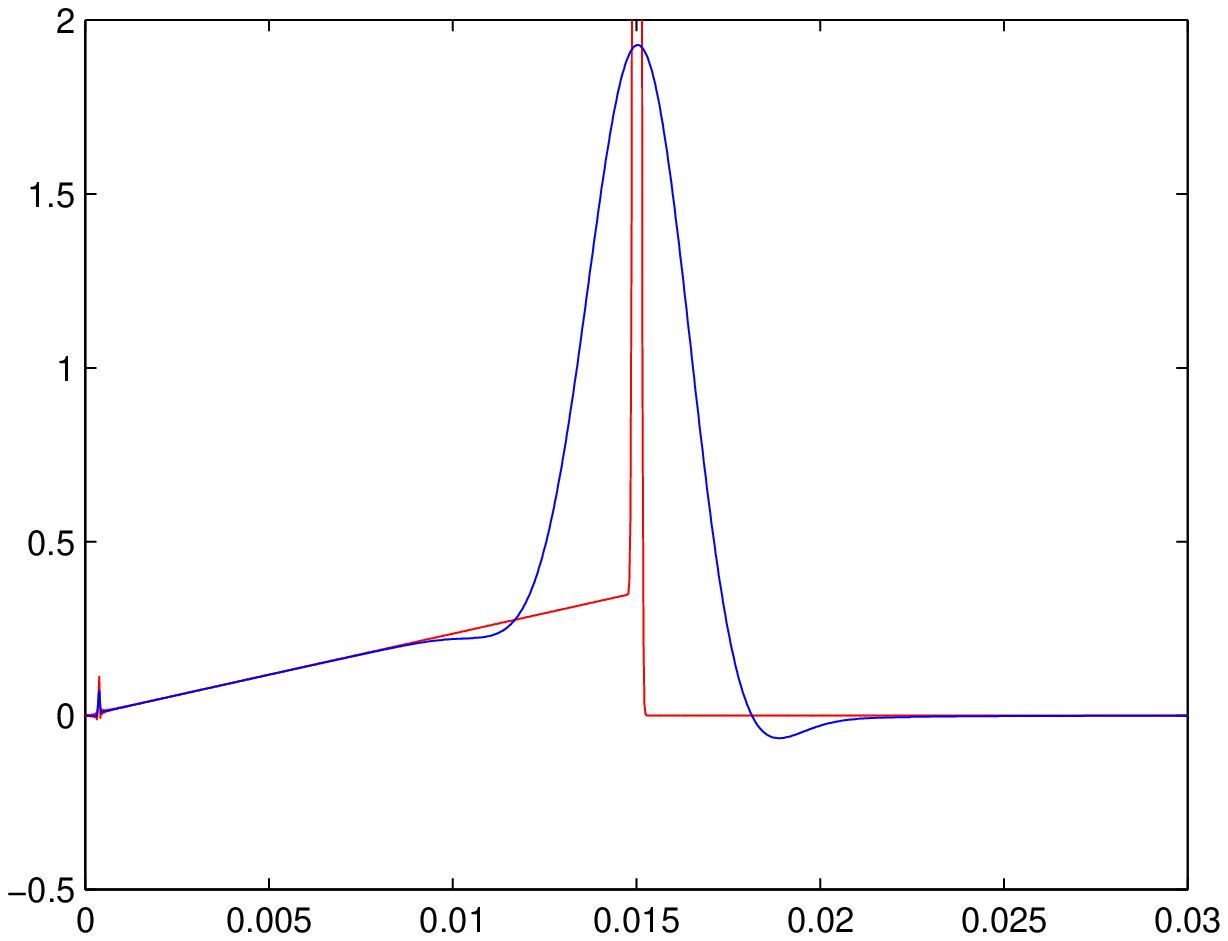}
\end{center}
\caption{Temporal response to a spatio-temporal delta function using a purely elastic Green's function (red line) and a 
viscous Green's function (blue line): Left, $y = 1.5$, $\nu_s = 4$ ; Center, $y = 2$, $\nu_s = 0.2$ ; Right, $y = 2.5$, $\nu_s = 0.002.$}
\end{figure}

 \begin{figure} 
\begin{center}  \label{fig:error_stationnary}
  \includegraphics[width=6cm]{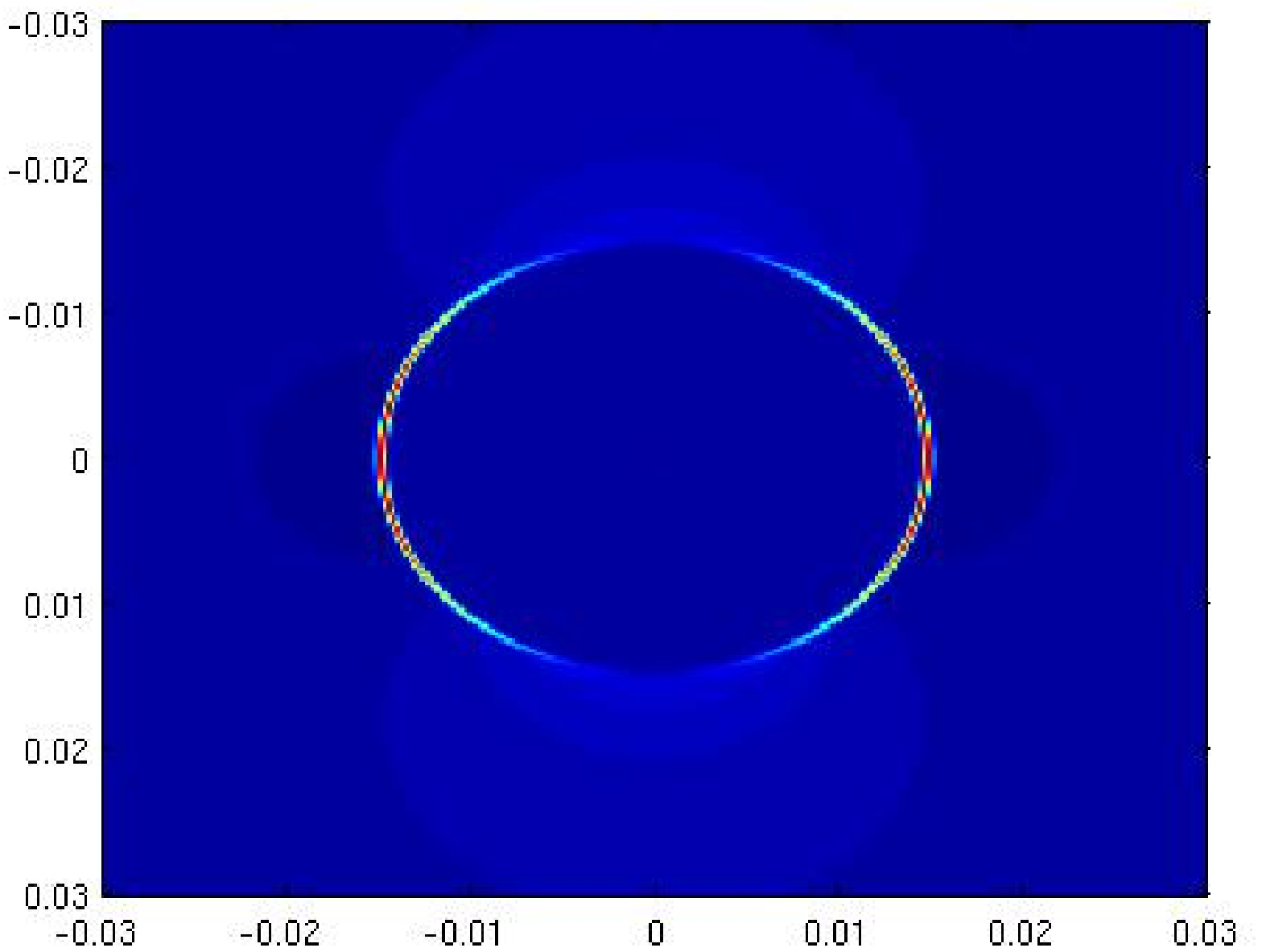}
 \includegraphics[width=6cm]{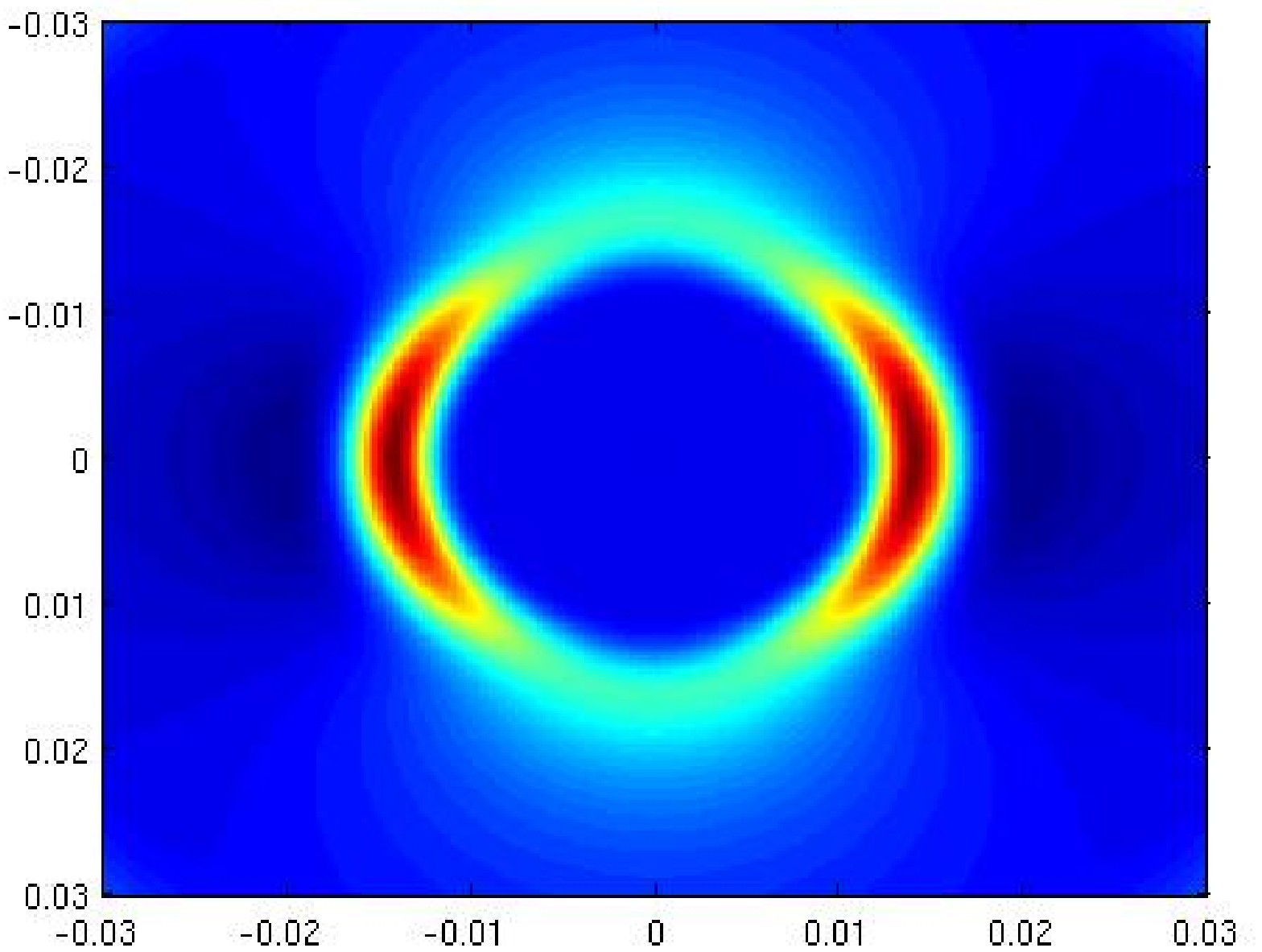}
 \includegraphics[width=6cm]{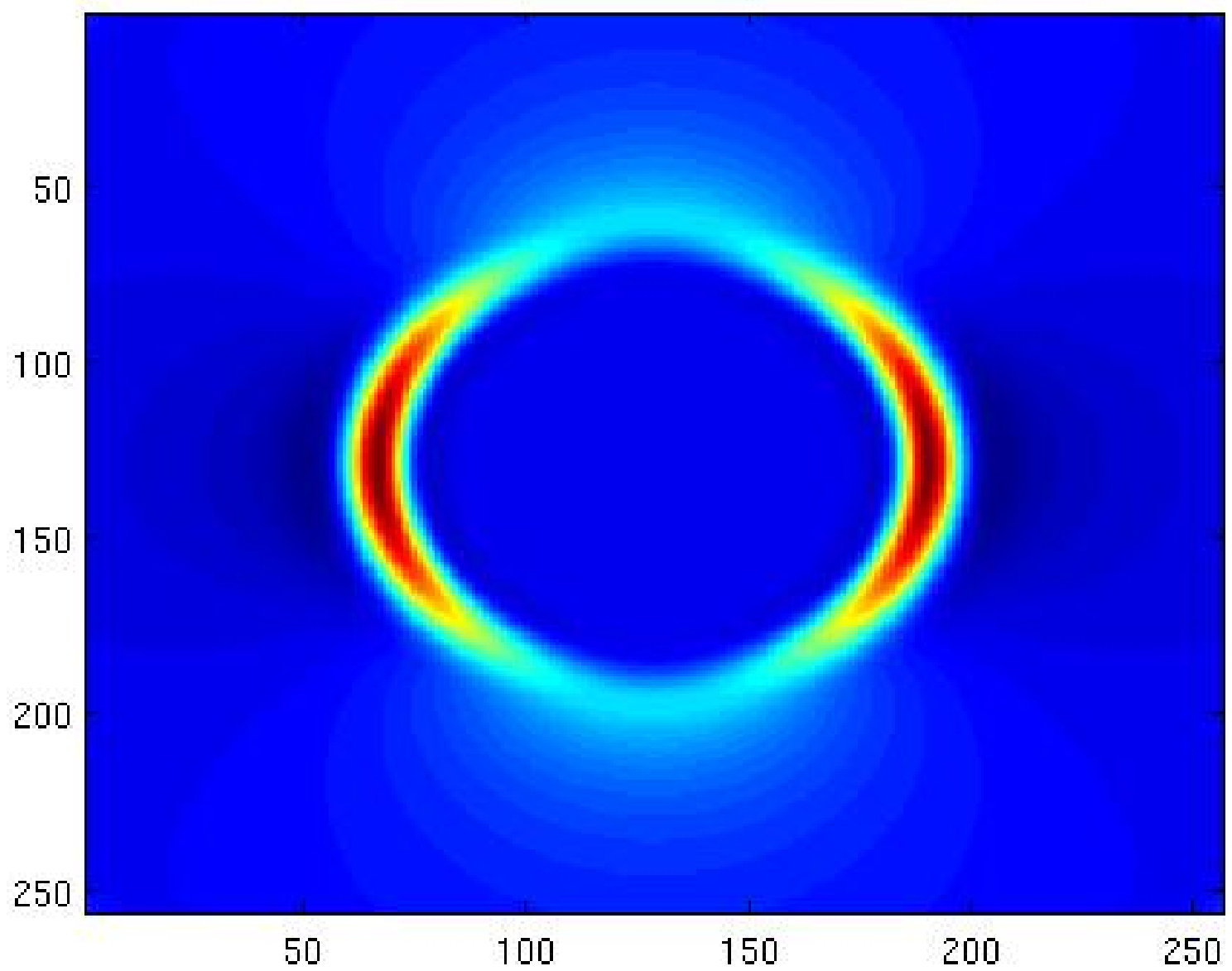}
\end{center}
\caption{ $2D$ spatial response to a spatio-temporal delta function  at $t =  0.015$ with a purely elastic Green's function,  a viscous Green's function  with $y = 2$, $\nu_s = 0.2$ and $y = 2.5$, $\nu_s = 0.002 $. }
\end{figure}

 \begin{figure} 
\begin{center}  \label{fig:phase_stationnaire}
  \includegraphics[width=6cm]{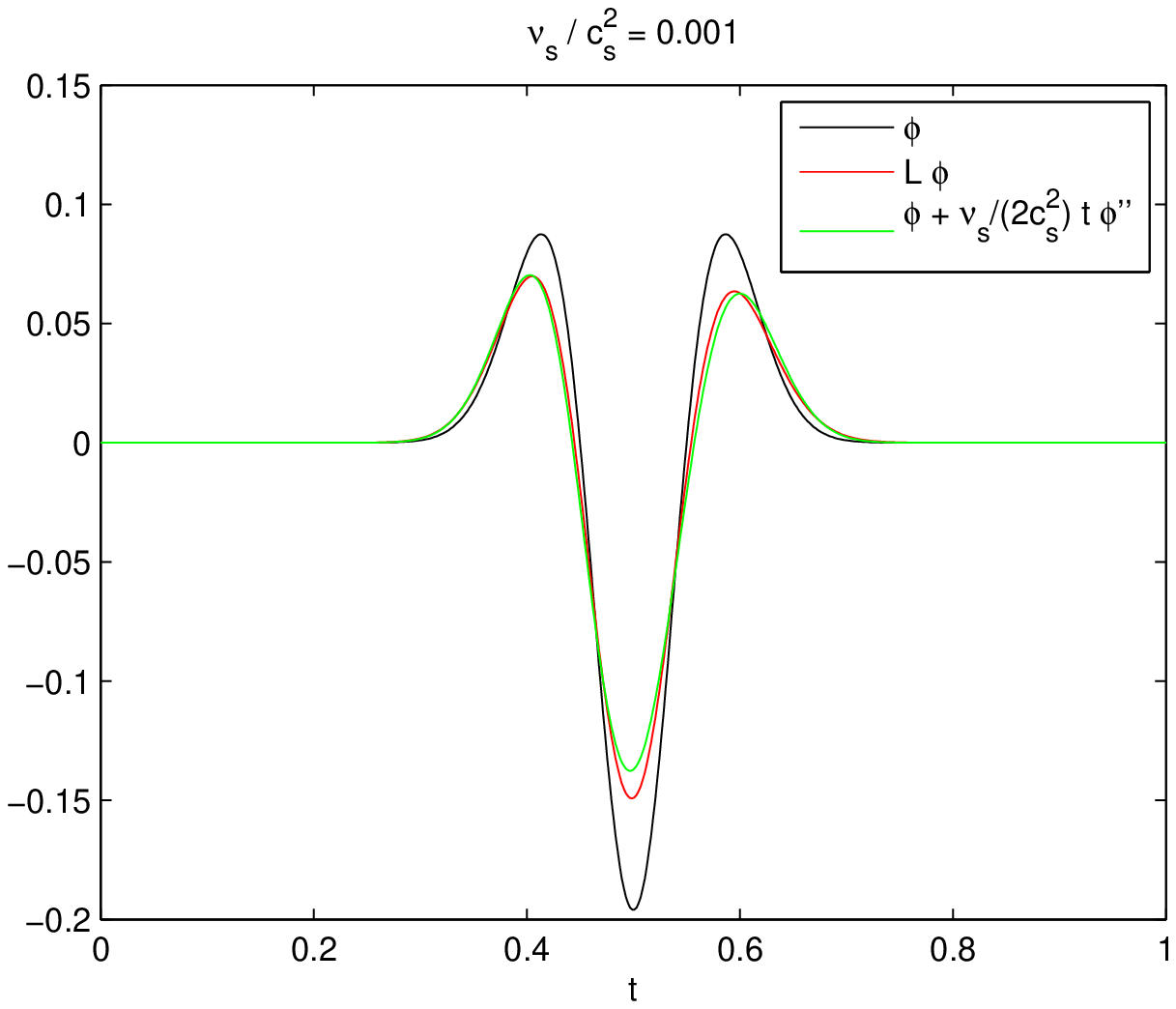}
 \includegraphics[width=6cm]{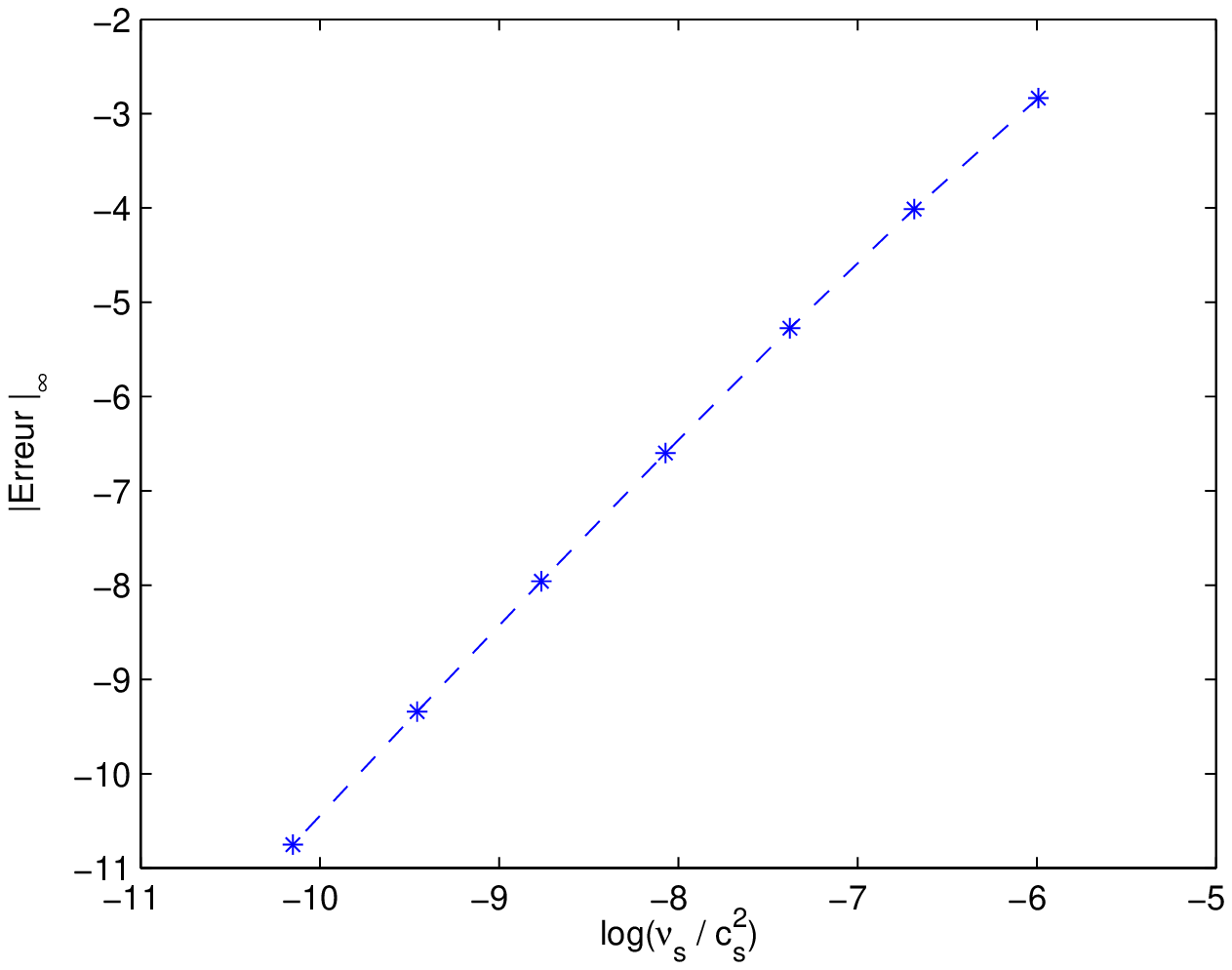}
\end{center}
\caption{Approximation of $L$ via stationary phase theorem : Left, comparison between $L \phi$ and  $\phi + \frac{\nu_s}{c_s^2 } t \phi''$ where $\frac{\nu_s}{c_s^2}=0.0001$ and $\phi$ is a smooth function. Right: error   $ \frac{\nu_s}{c_s^2 }  \to \|L \phi - \phi + \frac{\nu_s}{c_s^2 } \|_{\infty}$ in logarithmic scale. }
\end{figure}

%%%
%%%
%%%
%%%%%%%%%%%%%%%%%%%%%%%
\section{Conclusion}%%%
%%%%%%%%%%%%%%%%%%%%%%%
%%%
%%%
%%%
In this paper, we have computed the Green function in a
visco-elastic medium obeying a frequency power-law. For the Voigt
model, which corresponds to a quadratic frequency loss, we have
used the stationary phase theorem \ref{StationaryPhase} to
reconstruct the ideal Green function from the viscous  one by
solving an ODE. Once the ideal Green function is reconstructed,
one can find its source $\xi$ using the algorithms in
\cite{HabibEmergBioMedImag, HabibModelBioMed, HabibLiliAcoustic,
Habib-Lili10} such as time reversal, back-propagation, and
Kirchhoff Imaging. For more general power-law media, one can
recover the ideal Green function from the viscous one by inverting
a fractional derivative operator. This would be the subject of a
forthcoming paper.

\appendix
\section{Proof of Theorem (\ref{Thm:Approx})}\label{ApenA}

The proof of theorem (\ref{Thm:Approx}) is based on the
following theorem (see \cite[Theorem 7.7.1]{hormander}).
%%%
%%%
%%%
%%%%%%%%%%%%%%%%%%%%%%%%%%%%%%%%%%%%%%%%%%%%%%%%%%%%%%%%%%%%%%%%
\begin{thm}{\it\bf (Stationary Phase)}\label{StationaryPhase}%%%
%%%%%%%%%%%%%%%%%%%%%%%%%%%%%%%%%%%%%%%%%%%%%%%%%%%%%%%%%%%%%%%%
%%%
%%%
%%%
Let $K \subset [0,\infty)$ be a compact set,  $X$ an open
neighborhood of  $K$ and $k$ a positive integer. If $\psi \in
C_0^{2k}(K)$, $f \in C^{3k+1}(X)$ and $Im(f) \geq 0$ in $X$,
$Im(f(t_0)) = 0$, $f'(t_0) = 0$, $f''(t_0) \neq 0$, $f' \neq 0$ in
$K\setminus \{t_0\}$ then for $\epsilon > 0 $
%%%
%%%
%%%
%%%
\begin{eqnarray*}
 \left| \int_K \psi(t) e^{i f(t)/\epsilon} dx  - e^{i f(t_0)/\epsilon} \left( \lambda f''(t_0)/2\pi i \right)^{-1/2} \sum_{j<k} \epsilon^{j} L_j \psi \right| \leq C \epsilon^{k} \sum_{\alpha \leq 2 k } \sup |\psi^{(\alpha)}(x)|.
\end{eqnarray*}
%%%
%%%
%%ù
Here $C$ is bounded when $f$ stays in a bounded set in
$C^{3k+1}(X)$ and $|t-t_0|/|f'(t)|$ has a uniform bound. With,
$$ g_{t_0}(t) = f(t) - f(t_0) - \frac{1}{2}f''(t_0)(t - t_0)^2,$$
which vanishes up to third order at $t_0$, we have
$$ L_j \psi = \sum_{\nu - \mu = j} \sum_{2 \nu \geq 3 \mu} i^{-j}  \frac{2^{-\nu}}{\nu! \mu !}
(-1)^{\nu} f''(t_0)^{-\nu} (g_{t_0}^{\mu} \psi)^{(2
\nu)}(t_0).~~~~\Box$$
\end{thm}
%%%
%%%
%%%
Note that $L_1$ can be expressed as the sum  $L_1 \psi = L^{1}_1
\psi + L^{2}_1 \psi + L^{3}_1 \psi$, where $L_1^j$ is
respectively associate  to the pair $(\nu_j,\mu_j) =
(1,0),(2,1),(3,2)$ and is identified to
%%%
%%%
%%%
$$ \begin{cases}
 L^{1}_1 \psi &= \frac{-1}{2 i } f''(t_0)^{-1} \psi^{(2)}(t_0), \\
 L^{2}_1 \psi &= \frac{1}{2^2 2! i } f''(t_0)^{-2} (g_{t_0}u)^{(4)}(t_0) = \frac{1}{8 i } f''(t_0)^{-2} \left(g_{t_0}^{(4)}(t_0) \psi(t_0) + 4 g_{t_0}^{(3)}(t_0) \psi'(t_0)  \right) , \\
 L^{3}_1 \psi &=  \frac{-1}{2^3 2!3! i} f''(t_0)^{-3} (g_{t_0}^2 \psi)^{(6)}(t_0) = \frac{-1}{2^3 2!3! i}
 f''(t_0)^{-3} (g_{t_0}^2)^{(6)}(t_0) \psi(t_0).
 \end{cases}
$$
%%%
%%%
%%%

Now we turn to the proof of formula (\ref{formasympfinal}). Let us
first consider the case of operator $L^{*}$. We have
%%%
%%%
%%%
$$ \tilde L^{*} \phi(t) =  \int_0^{+\infty} \frac{\tau}{t}\phi(\tau)  \frac{c_s}{\sqrt{2 \pi \nu_s t}} e^{-\frac{c_s^2(\tau-t)^2}{2\nu_s t}}  d\tau = \frac{1}{t\sqrt{\epsilon}}\left(  \int_0^{+\infty}  \psi(\tau) e^{ i f(\tau)/\epsilon}\right),$$
%%%
%%%
%%%
with, $f(\tau) = i  \pi (\tau-t)^2$,  $\epsilon = \frac{ 2 \pi \nu_s t }{c_s^2}$ and $\psi(\tau) = \tau \phi(\tau)$.
Remark that the phase $f$ satisfies at $\tau = t$ ,  $f(t) = 0$, $f'(t) = 0$, $f''(t) =   2 i \pi \neq 0$.
Moreover, we  have
%%%
%%%
%%%
$$
\begin{cases}
 e^{i f(t)/\epsilon} \left( \epsilon^{-1} f''(t)/2 i \pi  \right)^{-1/2} =  \sqrt{\epsilon}   \\
 g_t(\tau) = f(\tau) - f(t) - \frac{1}{2}f''(t)(\tau - t)^2 =  0 \\
 L_1 \psi(t) = L^{1}_1 \psi(t) = \frac{-1}{2 i} f''(t)^{-1} \psi^{''}(t) = \frac{1}{4 \pi} (t \phi)''.
\end{cases}
$$
%%%
%%%
%%%
Thus, Theorem \ref{StationaryPhase} implies that
%%%
%%%
%%%
$$ \left| \tilde L^{*} \phi(t) - \left( \phi(t) +  \frac{\nu_s}{2 c_s^2 } (t \phi)''\right)  \right| \leq \frac{C}{t} \epsilon^{3/2} \sum_{\alpha \leq 4} \sup |(t \phi)^{(\alpha)}|.$$
%%%
%%%
%%%
The case of the operator $\tilde L$ is very similar. Note that
%%%
%%%
%%%
$$\tilde L \phi(t) = \int_0^{+\infty}  \frac{t}{\tau} \phi(\tau)  \frac{ c_s}{\sqrt{2 \pi\nu_s  \tau}} e^{-
 \frac{c_s^2(\tau-t)^2}{2\nu_s \tau}} \, d\tau = \frac{t}{\sqrt{\epsilon}} \left(  \int_0^{+\infty}  \psi(\tau) e^{ i f(\tau)/\epsilon}\right), $$
%%%
%%%
%%%
with $f(\tau) = i \pi \frac{(\tau-t)^2}{\tau}$,  $\epsilon = \frac{\nu_s}{2 \pi c_s^2}$ and $\psi(\tau) = \phi(\tau) \tau^{-\frac{3}{2}}$.
It follows that
%%%
%%%
%%%
$$f'(\tau) = i \pi \left( 1 - \frac{t^2}{\tau^2} \right), \quad f''(\tau) =  2 i \pi  \frac{t^2}{\tau^3}, \quad f''(t) =  2 i \pi  \frac{1}{t},$$
%%%
%%%
%%%
and the function $g_t(\tau)$  equals to
%%%
%%%
%%%
$$ g_{t}(\tau) =  i \pi \frac{(\tau - t)^2}{ \tau }  - i \pi \frac{(\tau-t)^2}{t} =  i\pi \frac{(t - \tau)^3}{ \tau t  }.$$
%%%
%%%
%%%
We deduce that
%%%
%%%
%%%
$$ \begin{cases}
     (g_{t} \psi )^{(4)}(t) &= \left( g_{t}^{(4)}(t) \psi(t) + 4 g_{t}^{(3)}(t) \psi'(t)  \right) =  i\pi \left(\frac{24}{t^3} \psi(t)-\frac{24}{t^2} \psi'(t)  \right)  \\
     (g_{t}^2 \psi )^{(6)}(t)&=  (g_{t}^2)^{(6)}(t) \psi(t) = - \pi^2 \frac{6 !}{t^4} \psi(t),
   \end{cases}
$$
%%%
%%%
%%%
and then,
%%%
%%%
%%%
$$ \begin{cases}
    L^1_1 \psi =   \frac{-1}{i} \left( \frac{1}{2}(f''(t))^{-1} \psi''(t) \right) = \frac{1}{4 \pi} t \left(\frac{\tilde{\phi}}{\sqrt{t}} \right)'' = \frac{1}{4\pi} \left( \sqrt{t} \tilde{\phi}''(t) - \frac{\tilde{\phi}'(t)}{\sqrt{t}} + \frac{3}{4} \frac{\tilde{\phi}}{t^{3/2}}   \right) \\
    L^2_1 \psi =   \frac{1}{8 i } f''(t)^{-2} \left(g_{t}^{(4)}(s) \psi(s) + 4 g_{t}^{(3)}(t) \psi'(t)  \right) = \frac{1}{4 \pi} \left( 3 \left( \frac{\tilde{\phi}(t)}{\sqrt{t}} \right)' - 3   \frac{\tilde{\phi}(t)}{t^{3/2}}   \right) = \frac{1}{4 \pi} \left( 3 \frac{\tilde{\phi}'(t)}{\sqrt{t}} - \frac{9}{2} \frac{\tilde{\phi}(t)}{t^{3/2}}  \right)\\
    L^3_1 \psi =   \frac{-1}{2^3 2!3! i} f''(t)^{-3} (g_{t}^2)^{(6)}(t) \psi(s) = \frac{1}{4 \pi} \left( \frac{15}{4}  \frac{\tilde{\phi}(t)}{t^{3/2}} \right),
\end{cases}
$$
%%%
%%%
%%%
where $\tilde{\phi}(\tau) =  \phi(\tau)/\tau$. Then, we have

%%%
%%%
%%%
\begin{eqnarray*}
 L^1 \psi &=&   L^1_1 \psi  +    L^2_1 \psi  + L^3_1 \psi  \\
       &=& \frac{1}{4\pi} \left(  \sqrt{t} \tilde{\phi}''(t) + \left( 3 - 1  \right) \frac{\tilde{\phi}'(t)}{\sqrt{t}} + \left(\frac{3}{4}  - \frac{9}{2}  + \frac{15}{4} \right)  \frac{\tilde{\phi}(t)}{t^{3/2}}  \right) = \frac{1}{4 \pi \sqrt{t}}  \left( t \tilde{\phi}(t) \right)'' =\frac{1}{4 \pi \sqrt{t}}  \phi''(t),
\end{eqnarray*}
%%%
%%%
%%%
and  again Theorem \ref{StationaryPhase} shows that
%%%
%%%
%%%
$$ \left| \tilde L \phi(t) - \left( \phi(t) +  \frac{\nu_s}{2 c_s^2} t \phi''(t)
 \right) \right|  \leq C t \epsilon^{3/2} \sum_{\alpha \leq 4} \sup |\psi^{(\alpha)}(t)|.$$
%%%
%%%
%%%
%%%%%%%%%%%%%%%%%%%%%%%%%%%%
{\it\bf\Large Acknowledgement\\\\}  {The authors would like to thank Prof. Habib Ammari for proposing this problem and for his fruitful pieces of advice. This work is supported by the foundation \emph{Digiteo} and the \emph{Higher Education Commission of Pakistan}}.% under the scholarship programe \emph{PhD in Natural and Basic Sciences-France}. }
%%%
%%%
%%%

\end{document}